\documentclass[12pt]{article}

\usepackage[latin1]{inputenc}
\usepackage{amsmath,amssymb}
\usepackage{amsmath,accents}
\usepackage{alltt,graphics,graphpap,color}
\usepackage[dvips]{graphicx}
\usepackage{amsfonts}
\usepackage{amscd}
\usepackage{mathrsfs}
\usepackage{hyperref}

\newtheorem{Theorem}{Theorem}[section]

\newtheorem{Proposition}[Theorem]{Proposition}
\newtheorem{Lemma}[Theorem]{Lemma}
\newtheorem{Corollary}[Theorem]{Corollary}
\newtheorem{Remark}[Theorem]{Remark}

\newtheorem{Notation}[Theorem]{Notation}

\newcommand{\rl}{{{\rm I} \kern -.15em {\rm R} }}
\newcommand{\nat}{{{\rm I} \kern -.15em {\rm N} }}
\newcommand{\sph}{{{\rm l} \kern -.45em {\rm S} }}

\newcommand{\mm}{{\cal M}}
\newcommand{\bmm}{{\overline {\cal M} }}

\newcommand{\nb}{\nabla}

\newcommand{\proof}{\emph{Proof. }}
\newcommand{\cvd}{\hfill$\square$ \bigskip}

\newcommand{\sss}{\mathbb S}
\newcommand{\qq}{\mathbb H}
\newcommand{\kk}{\mathbb K}
\newcommand{\pp}{\mathbb P}
\newcommand{\cc}{\mathbb C}

\newcommand{\tm}{T_{max}}

\newcommand{\aaa}{\left| A\right|^2}

\newcommand{\hhh}{\left| H\right|^2}

\newcommand{\bnb}{\overline\nb}

\newcommand{\dt}{\displaystyle{\frac{\partial}{\partial t}}}

\newcommand{\vc}{\mathscr V}

\newcommand{\hc}{\mathscr H}
\newcommand{\hl}{^{\mathscr H}}
\newcommand{\tc}{\mathcal T}
\newcommand{\ac}{\mathcal A}
\newcommand{\aaap}{\left|A'\right|^2}
\newcommand{\hhhp}{\left|H'\right|^2}
\newcommand{\tl}{^{\mathscr T}}
\newcommand{\bb}{{\cal B}}
\newcommand{\bbb}{{\overline {\cal B} }}

\hypersetup{linktoc=page,colorlinks=true,linkcolor=red,  citecolor=blue, }

\begin{document}
\def\qed{\hbox{\hskip 6pt\vrule width6pt height7pt
depth1pt  \hskip1pt}\bigskip}


\title{Mean curvature flow and Riemannian submersions}

\author{\sc G. Pipoli}
\date{}

\maketitle


{\small \noindent {\bf Abstract:} We give a sufficient condition ensuring that the mean curvature flow commutes with a Riemannian submersion and we use this result to create new examples of evolution by mean curvature flow. In particular we consider evolution of pinched submanifolds of the sphere, of the complex projective space, of the Heisenberg group and of the tangent sphere bundle equipped with the Sasaki metric.} \medskip

\noindent {\bf MSC 2010 subject classification} 53C44, 35B40 \bigskip

\section{Introduction}\hspace{5 mm}  \setcounter{equation}{0}\setcounter{Theorem}{0}\noindent

Let $F_0:\mm \to \left(\bmm,\bar g\right)$ be a smooth immersion of a $m$-dimensional manifold into a Riemannian manifold $\bmm$ of dimension $m+k$, called \emph{ambient space}. We denote by $A$ the second fundamental form and  by $H$ the mean curvature vector associated with the immersion. 
The evolution of $\mm_0 = F_0(\mm)$ by mean curvature  flow is the one--parameter family of  immersions $F:\mm \times [0,\tm[ \,\to\bmm$ satisfying
 \begin{equation}\label{1.1}
\left\{\begin{array}{l}
\dt F(p,t)= H,
\qquad p \in \mm, \, t \geq 0, \medskip \\
F(\cdot,0)=F_0.
\end{array}\right.
\end{equation}
It is well known that if $\mm_0$ is closed, then this problem has a uniquely defined smooth solution up to some maximal time $\tm \leq \infty$. Very often we identify the immersion $F(.,t)$ with the immersed submanifold $\mm_t=F(\mm,t)$.

Let $(\mm,g_{\mm})$ and $(\bb,g_{\bb})$ two Riemannian manifolds of dimension $m$ and $b$ respectively. A \emph{Riemannian submersion} is a smooth map $\pi:\mm\to \bb$ satisfying the following axioms $S1$ and $S2$.
\begin{itemize}
\item [S1)]$\pi$ has maximal rank;
\end{itemize}
For every $p\in \mm$, $\pi^{-1}(p)$ is a submanifold of $\mm$ called \emph{fiber} over $p$. A vector field on $\mm$ is called \emph{vertical} if it is always tangent to fibers, \emph{horizontal} if always orthogonal to fibers. The second axiom is
\begin{itemize}
\item [S2)]for every $X$, $Y$ horizontal vectors we have 
$$
g_{\mm}(X,Y)=g_{\bb}(\pi_* X,\pi_* Y)\circ \pi.
$$
\end{itemize}
$\mm$ is called \emph{total space} of the submersion and $\bb$ is called \emph{base}. Axiom $S1$ implies that $m\geq b$ and so the dimension of the fibers is $\hat m=m-b$. Axiom $S2$ says that $\pi$ preserves lengths of horizontal vectors.

The firs theorem proved in this paper explore the symmetries of the mean curvature flow and gives a sufficient condition ensuring that this flow commutes with a submersion. We consider submersions defined by the action of a group of isometries. Let $G$ be a Lie group acting as isometries of a Riemannian manifold $(\bmm, g_{\bmm})$. Suppose that the quotient space, obtained identifying the points of a orbit of the action of $G$, is a smooth manifold $\bbb=\bmm/G$ and consider the induced metric $g_{\bbb}$ on it. The natural projection $\pi:\bmm\to\bbb$ is a Riemannian submersion with fibers the orbits of $G$. If the action of $G$ is free we have the well-known principal bundles. In this case the fibers of $\pi$ are isometric to the group $G$. The best known examples of such submersions are probably the Hopf fibrations: $\pi_1:\sss^{2n+1}\to\cc\pp^n$ and $\pi_2:\sss^{4n+3}\to\qq\pp^n$.

 Lifting a submanifold of $\bbb$ we have a $G$-invariant submanifold of $\bmm$, vice versa projecting a $G$-invariant submanifold of $\bmm$ we get a submanifold of $\bbb$. We want to study how the mean curvature flow is related to a submersion.

\begin{Theorem}\label{mainsub}
Let $\pi:\bmm\to\bbb=\bmm /G$ be a Riemannian submersion. If $\pi$ has closed and minimal fibers then the mean curvature flow of any closed submanifold commutes with the submersion. More precisely let $\mm_0$ is a $G$-invariant submanifold of $\bmm$ and $\bb_0=\pi(\mm_0)$ then the mean curvature flow of $\mm_0$ and $\bb_0$ are defined up to the same maximal time $\tm$ and $\pi(\mm_t)=\bb_t$ for any time $0\leq t<\tm$. 

\end{Theorem}

Note that closedness of fibers and of the initial immersions guarantees the uniqueness of the solution of mean curvature flow of the submanifold $\bb_0$ and its lift. The mean curvature flow in manifolds with symmetries was studied by several authors, for example Pacini in \cite{Pa} considered the evolution of the orbits of a group of isometries. The proof of Theorem 1.1 is based on the fundamental equations for submersions which are derived in the classical paper by O'Neill \cite{O}. Although similar computations already appear in the previous literature on geometric flows, see in particular the paper \cite{S0}, it seems to us that this result was never explicitly observed before. The main part of this paper is devoted to the applications to specific examples, where we obtain new convergence results for the mean curvature flow by lifting to the ambient space the known theorems for the base manifold.

For example consider the main theorems of \cite{PS}: they concern evolution by mean curvature flow of pinched submanifolds of $\cc\pp^n$ and pinched hypersurfaces of $\qq\pp^n$. Lifting this result with Theorem \ref{mainsub} applied to the Hopf fibrations we have the following new examples of evolution of pinched hypersurfaces of the sphere.

\begin{Proposition}\label{hopf01=1}
Let $\mm_0$ be a closed $\sss^1$-invariant hypersurface of $\sss^{2n+1}(c)$, the sphere of constant sectional curvature $c>0$ with $n\geq 3$. If $\mm_0$ satisfies
\begin{equation}\label{eq_nguyen1}
\aaa < \frac{1}{2n-2}\hhh+4c,
\end{equation}
then the mean curvature flow of $\mm_0$ develops a singularity in finite time and converges to a $\sss^1$, therefore such an $\mm_0$ is diffeomorphic to a $\sss^1\times\sss^{2n-1}$.
\end{Proposition}


The pinching inequality \eqref{eq_nguyen1} is weaker than the one found by Huisken in \cite{H3}, but we have the further assumption about the $\sss^1$-invariance. In fact we do not find any of the two possibilities described by Huisken, i.e. the convergence to a round point in finite time or the convergence to a totally geodesic submanifold in infinite time. Moreover the pinching condition \eqref{eq_nguyen1} is the same studied by Nguyen in \cite{Ng}. We reached a similar result: a cylindrical singularity, 
but in our case we have a more complete result, in fact we found the global behavior of the evolution and not only around a singularity.

Another result of the same kind is the following.
\begin{Proposition}\label{hopf02=1}
Let $\mm_0$ be a closed $\sss^3$-invariant hypersurface of $\sss^{4n+3}(c)$, with $n\geq 3$. If $\mm_0$ satisfies
$$
\aaa < \frac{1}{4n-2}\hhh+8c,
$$
then the mean curvature flow of $\mm_0$ develops a singularity in finite time and converges to a $\sss^3$, then such an $\mm_0$ is diffeomorphic to a $\sss^3\times\sss^{4n-1}$.
\end{Proposition}

Note that $\sss^1$ is a subgroup of $\sss^3$, then if a submanifold of $\sss^{4n+3}$ is $\sss^3$-invariant, we can project it both to $\cc\pp^{2n+1}$ and $\qq\pp^n$. Putting together Propositions \ref{hopf01=1} and \ref{hopf02=1} we have a negative result.

\begin{Corollary}\label{hopf03=1}
There are no closed $\sss^{3}$-invariant hypersurfaces of $\sss^{4n+3}(c)$ such that 
$$
\aaa<\frac{1}{4n}\hhh+4c.
$$
\end{Corollary}

The paper is organized as follows. In section 2 we recall some notation and preliminary results, in particular we compute the relationships between the second fundamental forms of a submanifold and its lift through a Riemannian submersion. In Section 3 first we show that the invariance of a submanifold respect to a group of isometries of the ambient manifold is preserved by the mean curvature flow, then Theorem \ref{mainsub} is proved. In Section 4 we have the applications of Theorem \ref{mainsub}. Propositions \ref{hopf01=1} and \ref{hopf02=1} just discussed are proved in a more general setting: we deform the metric of the sphere using the canonical variations of the Hopf fibrations. After that some other examples are described. Proposition \ref{codalta} concerns $\sss^1$-invariant pinched submanifolds of higher codimension of the sphere: we found the alternative between the convergence in finite time to a $\sss^1$ and the convergence in infinite time to a totally geodesic submanifold, that is a $\sss^m$ for some $m$. This result extends the one of Baker \cite{Ba}. Proposition \ref{hopf04} is about pinched hypersurfaces of the complex projective space. The ambient space is the Heisenberg group in Propositions \ref{hei01} and \ref{hei02}. Finally in Proposition \ref{sasaki01} we study submanifolds of the tangent sphere bundle of the round spere equipped with the Sasaki metric and prove the alternative between the convergence in finite time to an orbit and the convergence in infinite time to a minimal, but not totally geodesic, limit.

\section{Preliminaries}\hspace{5 mm}  \setcounter{equation}{0}\setcounter{Theorem}{0}\noindent
In this section we recall some basic notions and fix some notations used through all this paper. Consider $F:\mm\rightarrow(\bmm, \bar g)$ a smooth immersion of an $m-$dimensional differential manifold $\mm$ into a Riemannian manifold $\bmm$ of dimension $m+k$. 
Unless told otherwise, geometric quantities of the submanifolds are indicated in the usual way, while for the ambient manifold we use a line over the common symbol. Moreover 
Latin letters $i,j,k,...$ are related to $T_x\mm$, the tangent space to $\mm$ at $x$, and Greek letters $\alpha, \beta, \gamma,...$ to the normal space $N_x\mm$. Fix $(x_1,\cdots,x_n)$ a local coordinate system around a point $x\in\mm$. The local expression of $g$ is
\begin{equation*}
\displaystyle{g_{ij}(x)=\bar g_{F(x)} \left(\frac{\partial F}{\partial x_i},\frac{\partial F}{\partial x_j}\right)}
\end{equation*}

\noindent Let  $\bnb$ be the Levi-Civita connection of $(\bmm,\bar g)$. The \emph{second fundamental form} $A$ of the immersion $F$ is defined for every $X$, $Y$ tangent vectors of $\mm$ by
$$
A(X,Y)=\left(\overline\nb_X Y\right)^{\perp},
$$
where $\perp$ denote the component normal to $\mm$.

Let $(\xi_1,\dots,\xi_k)$ be an orthonormal frame of $N_x\mm$, the second fundamental form can be written 
$$
A=h^{\alpha}\otimes \xi_{\alpha},
$$
where the $h^{\alpha}=\left(h^{\alpha}_{ij}\right)$ are symmetric $(0,2)$-tensors. Here and in the following, if there are no explicit signs of sum, we use Einstein notation, that is we sum over repeated indices.
The metric induces a natural isomorphism between tangent and cotangent space. In coordinates, this is expressed in terms of raising/lowering indexes by means of the matrices $g_{ij}$ and
$g^{ij}$, where $g^{ij}$ is the inverse of $g_{ij}$. The scalar product on the tangent space extends to any tensor bundle, by contracting any pair of lower and upper indices with $g_{ij}$ and $g^{ij}$ respectively. This also allows to define the norm of any tensor $T$. 
A function that we use very often is the norm of the second fundamental form
$$
\aaa=\sum_{\alpha}\left|h^{\alpha}\right|^2
$$
The trace respect to the metric $g$ of the second fundamental form is the \emph{mean curvature vector} $H$ :
$$
H=tr A=tr h^{\alpha}e_{\alpha}=g^{ij}h^{\alpha}_{ij}e_{\alpha}.
$$

\noindent It is independent of the orientation and it is well defined globally even if $\mm$ is non-orientable. Note that some authors defines the mean curvature as the trace of $A$ over $m$, of course this makes no substantial difference in the analysis.

 What follows is taken from the classical O'Neill's paper \cite{O}, many other interesting results about Riemannian submersions can be found in chapter 9 of \cite{Be2} and in the extensive monograph \cite{FIP}. 
Let $\pi:(\mm,g_{\mm})\to(\bb,g_{\bb})$ be a Riemannian submersion. 
If not specified otherwise, we use the same symbols for geometric quantities of $\mm$ and $\bb$. It will be clear from the context in which manifold we are. The same quantities of the fibers are distinguished by the superscript $\hat{\phantom{a}}$.
The \emph{vertical distribution} $\vc$ is the distribution of vertical vector fields, that is $\vc=\ker \pi_*$. Its orthogonal complement respect to $g_{\mm}$ is the \emph{horizontal distribution} $\hc$.  We denote with the same symbols $\hc$ and $\vc$ the projections of the tangent space of $\mm$ to the subspaces of horizontal and vertical vectors, respectively. Then every $X$ tangent to $\mm$ can be decomposed in an unique way in the sum of a horizontal and a vertical vectors:
$$
X=\hc X+\vc X.
$$
An horizontal vector field $X'$ is called \emph{basic} if there exists a vector fields $X$ on $\bb$ such that $\pi_*X'=X$,  in this case $X$ and $X'$ are said to be $\pi$-related.
There is an one-to-one correspondence between basic vector fields on $\mm$ and arbitrary vector fields on $\bb$: every basic vector field gives a vector field on $\bb$ by definition, while every $X$ tangent to $\bb$ has an unique horizontal lift $X\hl$ to $\mm$ characterized by $\pi_* X\hl=X$. Submersions are ruled by two tensors. For every $X$ and $Y$ tangent to $\mm$ we define
\begin{eqnarray*} 
\tc_{X}Y & = & \hc\nb _{\vc X}(\vc Y)+\vc\nb _{\vc X}(\hc Y);\\
\ac_{X}{Y} & = & \vc\nb_{\hc X}(\hc Y)+\hc\nb _{\hc X}(\vc Y).
\end{eqnarray*}
Note that if $X$ and $Y$ are tangent to fibers, i.e. vertical, then $\tc_{X}Y=\hat A(X,Y)$ the second fundamental form of the fibers as submanifolds of $\mm$. We have that $\tc\equiv 0$ if and only if each fiber is totally geodesic, while $\ac\equiv 0$ if and only if $\hc$ is integrable. 

Since we deal with the mean curvature flow we want to understand how a submanifold of $\bbb$ is related to its lift to $\bmm$: let $\pi: (\bmm, g_{\bmm})\to (\bbb, g_{\bbb})$ a Riemannian submersion, and $F:\bb\to\bbb$ an immersion. $\pi^{-1}(F(\bb))$ is a submanifold of $\bmm$ of the same codimension of $F(\bbb)$. Formally there is a manifold $\mm$, an immersion $F':\mm\to\bmm$ and a submersion that we indicate again with $\pi$, such that the following diagrams commutes:
$$
\begin{array}{ccc}
\bmm & \stackrel{\pi}{\longrightarrow} & \bbb \\
F' \uparrow & & \uparrow F\\
\mm& \stackrel{\pi}{\longrightarrow} & \bb
\end{array}
$$

We want to understand the link between $A$, the second fundamental form of $F$, and $A'$, the second fundamental form of $F'$. The main tool is the following O'Neill's formulas. 
\begin{Lemma}\label{o'neill1}\cite{O}
For every tangent vector fields on $\bbb$ $X$ and $Y$ we have
\begin{itemize}
\item [1)] $\left[X,Y\right]\hl=\hc\left[X\hl,Y\hl\right]$;
\item [2)] $\left(\bnb_X Y\right)\hl=\hc\left(\bnb_{X\hl}Y\hl\right)$.
\end{itemize}
\end{Lemma}

\begin{Lemma}\label{o'neill2}\cite{O}
Let $X$ and $Y$ be horizontal vector fields and $V$ and $W$ vertical vector fields. Then
\begin{itemize}
\item [1)] $ \bnb_{V}W=\tc_{V}W+\hat{\nb}_{V}{W}$;
\item [2)] $\bnb_{V}X=\hc\bnb_{V}X+\tc_{V}{X}$;
\item [3)] $\bnb_{X}V=\ac_{X}V+\vc\bnb_{X}{V}$;
\item [4)] $\bnb_{X}Y=\hc\bnb_{X}Y+\ac_{X}{Y}$.
\end{itemize}
\end{Lemma}

Note that, by construction, $\mm\equiv F'(\mm)$ is tangent to the fibers, then any vector normal to $\mm$ is necessarily horizontal. From Lemma \ref{o'neill1} and Gauss equation we have that for any $X$ and $Y$ tangent to $\mm$
\begin{equation}\label{aaa}
\begin{array}{rcl}
\bnb_{X\hl}Y\hl&=&\hc\left(\bnb_{X\hl}Y\hl\right)+\vc\left(\bnb_{X\hl}Y\hl\right)\\
 &=&\left(\bnb_X Y\right)\hl+\vc\left(\bnb_{X\hl}Y\hl\right)\\
 & = & \left(\nb_X Y\right)\hl+\left(A(X,Y)\right)\hl+\vc\left(\bnb_{X\hl}Y\hl\right).
\end{array}\end{equation}

By definition $A'(X\hl,Y\hl)=\left(\bnb_{X\hl}Y\hl\right)^{\perp}$ (the component normal to $\mm$), then it is an horizontal vector field. By \eqref{aaa} we have
$$
A'(X\hl,Y\hl)=\left(\left(\nb_X Y\right)\hl\right)^{\perp}+\left(\left(A(X,Y)\right)\hl\right)^{\perp}.
$$
The vector field $\left(\nb_X Y\right)\hl$ is the lift of a vector field tangent to $\bb$, then it is tangent to $\mm$. In the same way $\left(A(X,Y)\right)\hl$ is normal to $\mm$. Hence we have

\begin{equation}\label{horhor}
A'(X\hl,Y\hl)=\left(A(X,Y)\right)\hl .
\end{equation}

Now consider two vertical vector fields $V$ and $W$. They are tangent to $\mm$ by construction. $A'(V,W)$ is normal to $\mm$ then it is a horizontal vector field. By Lemma \ref{o'neill2} we have
\begin{equation}\label{verver}
A'(V,W)=\left(\bnb _{V}W\right)^{\perp}=\left(\hc\bnb _{V}W\right)^{\perp}=(\tc_{V}W)^{\perp}=(\hat A(V,W))^{\perp}.
\end{equation}
 Lemma \ref{o'neill2} does not say anything about the mixed terms $A'(X,V)$ with $X$ horizontal and $V$ vertical: they strongly depend on the specific submersion considered as we will see in the examples of section 4.

\begin{Notation}
For any submersion $\pi$ considered below $(X_1,\dots,X_m)$ denote a local orthonormal frame tangent to a submanifold of the base space around a point $p$ and $(V_1,\dots,V_{\hat m})$ is a local orthonormal set of vertical vector fields. Then around any point $q$ of the fiber $\pi^{-1}(p)$ we use the orthonormal basis $(X_1\hl,\dots,X_m\hl,V_1,\dots V_{\hat m})$ tangent to the lift of the submanifold. Moreover $(\xi_1,\dots,\xi_k)$ is a local orthonormal frame normal to a submanifold of the base, then $(\xi_1\hl,\dots,\xi_k\hl)$ is a local orthonormal frame normal to the lift of the submanifold considered. 
\end{Notation}\label{notazione_esempi}

\noindent Summarizing what we found we have
\begin{equation}\label{schemaA'}
A'= \left(\begin{array}{c|c}
h_{ij}\hl & \text{mixed terms}\\
\hline
\text{mixed terms} & \hat h_{ij}^{\perp}
\end{array}\right)
\end{equation}  
where $h_{ij}=A(X_i,X_j)$ and $\hat h_{ij}=\hat A (V_i,V_j)$. 

Starting from a fixed Riemannian submersion $\pi:(\mm,g_{\mm})\to(\bb,g_{\bb})$ there is a standard way to deform the metric $g_{\mm}$ to obtain again a Riemannian submersion. The \emph{canonical variation} of $g_{\mm}$ is the family of metrics $\left\{g_{\lambda}\right\}_{\lambda>0}$ on $\mm$ such that
$$
\begin{array}{lll}
g_{\lambda}(U,V)=\lambda g_{\mm}(U,V) &\text{if} & U,V\in\vc,\\
g_{\lambda}(X,Y)= g_{\mm}(U,V) &\text{if} & X,Y\in\hc,\\
g_{\lambda}(U,X)=0 &\text{if} & U\in\vc,\ X\in\hc.
\end{array}
$$
Obviously $g_1=g_{\mm}$. For any $\lambda>0$, $g_{\lambda}$ makes $\pi$ a Riemannian submersion with the same horizontal and vertical distributions and the same fibers. Let $\nb^{\lambda}$ be the Levi-Civita connection of the metric $g_{\lambda}$. A straightforward computation gives:
\begin{equation}\label{nblambda}
\begin{array}{l}
\vc\left(\nb^{\lambda}_UV\right)=\vc\left(\nb^1_UV\right),\quad \hc\left(\nb^{\lambda}_UV\right)=\lambda\hc\left(\nb^1_UV\right),\\
\nb^{\lambda}_XU=\nb^1_XU,\quad\nb^{\lambda}_UX=\nb^1_UX,\quad\nb^{\lambda}_XY=\nb^1_XY,
\end{array}
\end{equation}
for every $U,V\in\vc$ and $X,Y\in\hc$. It follows that $\pi:(\mm,g_{\mm}=g_1)\to(\bb,g_{\bb})$ has minimal (resp. totally geodesic) fibers if and only if $\pi:(\mm,g_{\lambda})\to(\bb,g_{\bb})$ has minimal (resp. totally geodesic) fibers for every $\lambda>0$. Moreover let $(V_1,\dots V_{\hat m})$ is a local $g_{\mm}$-orthonormal set of vertical vectors, then for any $\lambda>0$ $(\lambda^{-\frac 12}V_1,\dots \lambda^{-\frac 12}V_{\hat m})$ is $g_{\lambda}$-orthonormal. Using \eqref{nblambda} is easy to see that, respect to this basis, the equation \eqref{schemaA'} becomes

\begin{equation}\label{schemaA'var}
A'_{\lambda}= \left(\begin{array}{c|c}
h_{ij}\hl & \lambda^{-\frac 12}\text{ mixed terms}\\
\hline
\lambda^{-\frac 12}\text{ mixed terms} & \hat h_{ij}^{\perp}
\end{array}\right).
\end{equation}

\section{Symmetries of the mean curvature flow}\setcounter{Theorem}{0} \setcounter{equation}{0}
In this section we prove Theorem \ref{mainsub}.

\begin{Lemma}\label{sub01}
Let $F_0:\mm\to\bmm$ be a closed immersion and $\varphi$ an isometry of $\bmm$, then $\varphi$ commutes with the mean curvature flow. Formally if $G_0=\varphi\circ F_0$ and $F_t$ and $G_t$ are the evolutions of $F_0$ and $G_0$ respectively, we have that $G_t=\varphi\circ F_t$ for any time $t$ that the flow is defined.
\end{Lemma}
\proof Since $\varphi$ is an isometry we have
$$
\dt (\varphi\circ F_t)(p,t)=\varphi_*H^F(p,t)=H^{\varphi\circ F}(p,t),
$$
where $H^{\psi}$ is the mean curvature vector of $\psi$ for any immersion $\psi$.
Then $\varphi\circ F_t$ is a solution of the mean curvature flow of initial data $\varphi\circ F_0=G_0$. For the uniqueness of the solution we have the thesis. \cvd

It follows immediately that

\begin{Corollary}\label{sub02}
Let $F_0$ and $\varphi$ like in the statement of Lemma \ref{sub01} and $G$ a group of isometries of $\bmm$. We have
\begin{itemize}
\item [1)] if $F_0$ is $\varphi$-invariant, then $F_t$ is $\varphi$-invariant for any $t$,
\item [2)] if $F_0$ is $G$-invariant then $F_t$ is $G$-invariant for every time $t$. 
\end {itemize}
\end{Corollary}

\noindent\emph{Proof of Theorem \ref{mainsub}.} 
Let $F_0:\bb\to\bbb$ and $F'_0:\mm\to\bmm$ two immersions for $\bb_0$ and $\mm_0$ respectively. By hypothesis we have that $F'_0$ is $G$-invariant and $\pi\circ F'_0=F_0\circ\pi$.
The crucial point is that, since the fibers are minimal, we have that $H'$ is basic and is $\pi$-related with $H$, where $H$ is the mean curvature vector of any submanifold of $\bbb$ and $H'$ is the mean curvature vector of its lift to $\bmm$. In fact $H'$ is horizontal because it is normal to $\mm$. Using the notation of \eqref{schemaA'} we have
\begin{eqnarray*}
H' =tr A'& = & \sum_i A'(X_i\hl,X_i\hl)+\sum_iA'(V_i,V_i)\\
 & = & \sum_i A(X_i,X_i)\hl+\sum_i\hat A(V_i,V_i)^{\perp}\\
 & = & \left(\sum_i A(X_i,X_i)\right)\hl+\left(\sum_i\hat A(V_i,V_i)\right)^{\perp}\\
 & = & H\hl+\hat H^{\perp}.
\end{eqnarray*}
If the fibers are minimal we get $H'=H\hl$, that is $H'$ and $H$ are $\pi$-related. In particular $\pi_* H'=H$ holds.
Now let $F_t$ the evolution of $F_0$, $F'_t$ the lift of $F_t$, $\widetilde F'_t$ the evolution of $F'_0$ and $\widetilde F_t$ the projection of $\widetilde F'_t$ and $H$, $H'$, $\widetilde H'$ and $\widetilde H$ the respective mean curvature vectors. By construction we have that for any $t$
\begin{equation}\label{sub03}
\pi\circ F'_t=F_t\circ\pi,
\end{equation}
and $F'_t$ is $G$-invariant. Then in particular $H'$ is horizontal. Differentiating \eqref{sub03} we have 
$$
\pi_*\dt F'_t = \dt (F_t\circ \pi)=H=\pi_* H'.
$$
Then $\dt F'_t =H'+V'$ for some vertical vector field $V'$. Since $F'_t$ is $G$-invariant, $V'$ is tangent to $F'_t(\mm')$. Therefore
$$
\left(\dt F'_t\right)^{\perp} =H'.
$$
This means that, up to a tangential diffeomorphism, $F'_t$ is the solution of the mean curvature flow of initial data $F'_0$. Then $F'_t(\mm)=\widetilde F'_t(\mm)$ for every time $t$. Vice versa 
$$
\dt \left(\widetilde F_t\circ\pi\right) = \dt\left(\pi\circ\widetilde F'_t\right)=\pi_*\dt\widetilde F'_t=\pi_*\widetilde H'.
$$
Corollary \ref{sub02} says that $\widetilde F'_t$ is $G$-invariant as its initial data $F'_0$, then $\pi_*\widetilde H'=\widetilde H$. Then $\widetilde F_t$ is the evolution of initial data $F_0$, that is $\widetilde F_t(\bb)=F_t(\bb)$ for any time $t$. \cvd

\begin{Remark}\label{non_unico}
If fibers are not closed we do not know if the solution of the mean curvature flow of the lift is unique, but if they are minimal, the same proof given for Theorem \ref{mainsub} shows that the lift of the mean curvature flow is, in any case, a $G$-invariant solution of the mean curvature flow. In the same way the projection of a $G$-invariant solution is again an evolution by mean curvature. Then if the projection of the initial data $\mm_0$ is a closed submanifold $\bb_0$ then there exists only one $G$-invariant solution of initial data $\mm_0$.
\end{Remark}

\section{Examples and applications}\setcounter{Theorem}{0} \setcounter{equation}{0}

One of the best known examples of submersions is the family of the Hopf fibrations. Let $\kk$ be one of the field $\cc$ or the associative algebra $\qq$ and $a$ be the real dimension of $\kk$. We denote with $\mathbb{S}^n(c)$ the $n$-dimensional sphere with the canonical metric of constant curvature $c>0$. The action $T:\sss^{a-1}(1)\times\sss^{na+a-1}(c)\rightarrow\sss^{na+a-1}(c)$, $(\lambda,z)\mapsto\lambda z$ is by isometries which acts transitively on the fiber. The Hopf fibrations are $\pi:\sss^{na+a-1}(c)\rightarrow\kk\pp^n\equiv\sss^{na+a-1}/\sss^{a-1}$, $z\mapsto\left[z\right]$, where $\left[z\right]$ is the class of $z$ under the action $T$. The Riemannian metric that we consider on $\kk\pp^n$ is the one induced from the metric of $\sss^{na+a-1}(c)$ such that $\pi$ becomes a Riemannian submersion. For $\kk=\cc$ it is the well-known Fubini-Study metric. 

Let us consider first the Hopf fibration $\pi:\sss^{2n+1}\to\cc\pp^n$.
In this case $V=J\nu$ is the vertical unit vector field, where $J$ is the complex structure of $\cc^{n+1}$ and $\nu$ is the outward normal unit vector field of the sphere as submanifold of $\mathbb R^{2n+2}\equiv \cc^{n+1}$. Let $\bb_0$ a submanifold of $\cc\pp^n$ of dimension $m$ and codimension $k$ and  $\mm_0$ its lift to $\sss^{2n+1}$. The fibers $\sss^1$ are geodesics, hence of course minimal. For every $i$ $J(\xi_i)\hl$ is tangent to the sphere and horizontal. Define $J(\xi_i)\hl=-U_i+N_i$ where $U_i$ is the component tangent to $\mm_0$, while $N_i$ is normal.  We want to explicit the mixed terms in \eqref{schemaA'} for this submersion. As shown in \cite{O} for every horizontal lift we have 
$$
\ac_{X\hl}V=
J\left(X\right)\hl,
$$
If $X$ is tangent to $\bb_0$ then $A'(X\hl,V)$ is an horizontal vector field then, by Lemma \ref{o'neill2} and some trivial computation
\begin{eqnarray*}
A'(X\hl,V) & =& 
 \sum_i\bar g(X\hl,U_i)\xi_i\hl.
\end{eqnarray*}
Moreover since the fibers are geodesic curves, $A'(V,V)=0$, together to \eqref{horhor} that holds in general we have
$$
\left|A'\right|^2=\aaa+2\sum_i \left|U_i\right|^2.
$$

The canonical variation of the Hopf fibration gives a family $\left\{\bar g_{\lambda}\right\}_{\lambda>0}$ of metrics on $\sss^{2n+1}$. Respect to this metric, a unit vertical vector field is $V_{\lambda}=\lambda^{-\frac 12}J\nu$, then with the same computation seen above we have that
$$
\left|A'\right|^2=\aaa+2\lambda^{-\frac 12}\sum_i \left|U_i\right|^2.
$$
Then, for any $\lambda>0$
$$
\aaa\leq\aaap\leq\aaa+2 \lambda^{-\frac 12}cod\mm_0=\aaa+2\lambda^{-\frac 12} cod\bb_0
$$
holds. 
Obviously since $H'$ and $H$ are $\pi$-related we have that $\hhhp=\hhh$ in any case. 

In the same way we can study Hopf fibration $\pi:\sss^{4n+3}\to\qq\pp^n$. The fibers are $\sss^3$ which are totally geodesic. Let $J_1$, $J_2$ and $J_3$ the complex structures of $\qq^{n+1}$ given by the multiplication of the quaternionic imaginary units. Then $(V_1=J_1\nu,V_2=J_2\nu,V_3=J_3\nu)$ is an orthonormal basis of $\vc$. 
Following the same notations and the same computations of the previous case we define for every $i$ and $\alpha$ define $J_{\alpha}\xi_i\hl=-U_{i\alpha}+N_{i\alpha}$ where $U_{i\alpha}$ is tangent to $\mm_0$, while $N_{i\alpha}$ is normal.  Moreover
$$
\ac_{X\hl}V_{\alpha}=J_{\alpha}\left(X\right)\hl,
$$
and for every $\alpha$ and $\beta$
\begin{eqnarray*}
A'(X\hl,V_{\alpha}) & =& \sum_i\bar g(X\hl,U_{i\alpha}\hl)\xi_i\hl,\\
A'(V_{\alpha},V_{\beta}) & = & 0.
\end{eqnarray*}
Then we get
$$
\left|A'\right|^2=\aaa+2\sum_{i,\alpha} \left|U_{i\alpha}\hl\right|^2=\aaa+2\sum_{i,\alpha}\left|U_{i\alpha}\right|^2.
$$
The canonical variation of this Hopf fibration gives a second family $\left\{\tilde g_{\lambda}\right\}_{\lambda>0}$ of metric on $\sss^{4n+3}$. Likewise to the previous case we have
$$
\left|A'\right|^2=\aaa+2\lambda^{-\frac 12}\sum_{i,\alpha}\left|U_{i\alpha}\right|^2.
$$

\noindent Then for every $\lambda>0$ 
$$
\aaa\leq\aaap\leq\aaa+6 \lambda^{-\frac 12}cod\mm_0=\aaa+6 \lambda^{-\frac 12}cod\bb_0.
$$

As application of Theorem \ref{mainsub} we have the following results: as particular case, for $\lambda=1$, we have Proposition \ref{hopf01=1}

\begin{Proposition}\label{hopf01}
Let $\mm_0$ be a closed $\sss^1$-invariant hypersurface of $(\sss^{2n+1},\bar g_{\lambda})$, with $n\geq 3$. If $\mm_0$ satisfies
\begin{equation}\label{eq_nguyen}
\aaap < \frac{1}{2n-2}\hhhp+2+2\lambda^{-\frac 12},
\end{equation}
then the mean curvature flow of $\mm_0$ develops a singularity in finite time and converges to a $\sss^1$, then such an $\mm_0$ is diffeomorphic to a $\sss^1\times\sss^{2n-1}$.
\end{Proposition}

\proof Since $\mm_0$ is $\sss^1$ invariant we can project it with the Hopf fibration to an hypersurface $\bb_0$ of $\cc\pp^n$. For hypersurfaces we have necessarily $\aaap=\aaa+2\lambda^{-\frac 12}$. Then $\bb_0$ satisfies
$$
\aaa < \frac{1}{2n-2}\hhh+2
$$ 
By Theorem 1.1 of \cite{PS} the evolution of $\bb_0$ converges in finite time to a round point $p$. By Theorem \ref{mainsub} we have  that the evolution of $\mm_0$ convergences in finite time to the lift of this point, that is a fiber. \phantom{aa}\hfill \cvd


For higher codimension we have the following result.

\begin{Proposition}\label{codalta}
Consider $\mm_0$ a closed $\sss^1$-invariant submanifold of $(\sss^{2n+1},\bar g_{\lambda})$ of dimension $m$ and codimension $2\leq k<\frac{2n-3}{5}$ satisfying the pinching condition
\begin{equation}\label{pinchingcodalt}
\aaa<\frac{1}{m-2}\hhh+\frac{m-4-4k}{m-1}.
\end{equation}
Then either
\begin{itemize}
\item [1)] the evolution of $\mm_0$ converges in finite time to a $\sss^1$,
\item[] or
\item [2)] the evolution of $\mm_0$ is defined for any time $0\leq t<\infty$ and converges to a smooth totally geodesic submanifold, that is an $\sss^{2n-k+1}$.
\end{itemize}
If $k$ is odd only case 1) can occur.
\end{Proposition}
\proof The proof is the same of the previous Proposition: it follows again from Theorem 1.1 of \cite{PS} and Theorem \ref{mainsub}: using inequality $\aaa\leq \aaap$ we have that $\bb_0=\pi(\mm_0)$ satisfies the same pinching inequality \eqref{pinchingcodalt} and the thesis follows since $$\pi^{-1}\left(\cc\pp^{n-\frac k2}\right)=\sss^{2n-k+1}$$ \cvd

This time let us consider the canonical deformation of the Hopf fibration $\pi:\sss^{4n+3}\to\qq\pp^n$, for $\lambda=1$ we have Proposition \ref{hopf02=1}. 
\begin{Proposition}\label{hopf02}
Let $\mm_0$ be a closed $\sss^3$-invariant hypersurface of $(\sss^{4n+3},\tilde g_{\lambda})$, with $n\geq 3$. If $\mm_0$ satisfies
$$
\aaap < \frac{1}{4n-2}\hhhp+2+6\lambda^{-\frac 12},
$$
then the mean curvature flow of $\mm_0$ develops a singularity in finite time and converges to a $\sss^3$, then such an $\mm_0$ is diffeomorphic to a $\sss^3\times\sss^{4n-1}$.
\end{Proposition}
\proof For hypersurfaces we have $\aaap=\aaa+6\lambda^{-\frac 12}$, then $\bb_0=\pi(\mm_0)$ satisfies $\aaa<\frac{1}{4n-2}\hhh+2$, then by Theorem 7.1 of \cite{PS} the evolution of $\bb_0$ shrinks to a round point in finite time. The thesis follows as in the previous Propositions. \cvd






A further example is given by the submersion $\rho:\cc\pp^{2n+1}\to\qq\pp^n$ described in \cite{E2}: it is the submersion that makes commutative the following diagrams
\begin{equation}\label{diagrammarho}
\begin{array}{rcl}
\sss^{4n+3} & \stackrel{\pi_1}{\longrightarrow} &\cc\pp^{2n+1}\\
\pi_2\searrow & & \swarrow\rho\\
& \qq\pp^n&
\end{array}
\end{equation}
where $\pi_1$ and $\pi_2$ are the usual Hopf fibrations. The fibers of $\rho$ are $\cc\pp^1\equiv\sss^2(4)$ and hence they are totally geodesic. Lifting an hypersurface of $\qq\pp^n$ to an hypersurface of $\cc\pp^{2n+1}$ via $\rho$ we have that 
$\aaap=\aaa+4$. In the same way of the previous propositions we can prove the following result.

\begin{Proposition}\label{hopf04}
Let $\mm_0$ be a closed $\cc\pp^1$-invariant hypersurface of $\cc\pp^{2n+1}$. If $\mm_0$ satisfies
$$
\aaa < \frac{1}{4n-2}\hhh+6,
$$
then the mean curvature flow of $\mm_0$ develops a singularity in finite time and converges to a fiber $\cc\pp^1$, then such an $\mm_0$ is diffeomorphic to a $\sss^2\times\sss^{4n-1}$.
\end{Proposition}


The examples seen before are all principal bundles with compact fibers. An interesting case with non-compact fibers comes from the Heisenberg group $\qq^n$ (not to be confused with the algebra of quaternions!). The Heisenberg group is the Lie group  
$\mathbb R^{2n}\times \mathbb R$ endowed with the following product:
$$
(x,y,z)(x',y',z')=\left(x+x', y+y', z+z'+\frac{1}{2}\left(\left\langle x,y'\right\rangle-\left\langle y,x'\right\rangle\right)\right),
$$
where $x,x',y,y'\in\mathbb R^n$, $z,z'\in\mathbb R$ and $\left\langle \cdot,\cdot\right\rangle$ is the Euclidean scalar product of $\mathbb R^n$. Respect to the coordinates $(x,y,z)=(x_1,\dots,x_n,y_1,\dots,y_n,z)$ we define the following left invariant vector fields on $\qq^n$:
$$
\begin{array}{l}
\displaystyle{X_j = \frac{\partial}{\partial x_j}-\frac{1}{2}y_j\frac{\partial}{\partial z},}\quad
\displaystyle{Y_j = \frac{\partial}{\partial y_j}+\frac{1}{2}x_j\frac{\partial}{\partial z},}\quad
\displaystyle{V=\frac{\partial}{\partial z}.}
\end{array}
$$
Declaring orthonormal the basis $\left(X_j,Y_j,V\right)_j$ we have a left-invariant metric $\bar g$ on $\qq^n$. On $\cc^n$ consider the Euclidean metric, then $$\pi:(x,y,z)\in\qq^n\mapsto (x+iy)\in\cc^n$$ is a Riemannian submersion. The fibers are the vertical line: $$\pi^{-1}(x_0+iy_0)=\left\{\left.\left(x_0,y_0,t\right)\right|t\in\mathbb R\right\}.$$ Moreover $\vc=span\left\langle V\right\rangle$ and $\hc=span\left\langle X_j,Y_j\right\rangle_{j=1,\dots,n}$. The structural group is the group of vertical translations, that is the multiplication by a point of the type $(0,0,t)$. It is a group of isometries and it is isomorphic to $(\mathbb R,+)$. The Levi-Civita connection associated to $\bar g'$ is determined by
$$
\begin{array}{rcccl}
\displaystyle{\bnb_{X_j}{Y_j}} & = &\displaystyle{ -\bnb_{Y_j}{X_j} }& =&\displaystyle{\frac{1}{2}V,} \\
\displaystyle{\bnb_{X_j}{V} }& = &\displaystyle{ \bnb_{V}{X_j} }& =&\displaystyle{-\frac{1}{2}Y_j }\\
\displaystyle{\bnb_{Y_j}{V} }& = & \displaystyle{\bnb_{V}{Y_j} }& =&\displaystyle{ \frac{1}{2}X_j}
\end{array}
$$
and zero for all others pairs of vector of the basis $\left(X_j,Y_j,V\right)_{j=1,\dots,n}$. A proof can be found in \cite{Ma}. In particular $\bnb_{V}V$ vanishes, hence the fiber of $\pi$ are geodesics. On the horizontal distribution $\hc$ we have a complex structure $J$ defined on the vector of the basis by $JX_j=Y_j$ and $JY_j=-X_j$ for all $j$. Then more succinctly, for any horizontal vector field $Z$ on $\qq^n$ we have
\begin{equation}\label{nbheis}
\bnb_{Z}V=\bnb_{V}Z=-JZ.
\end{equation}
Now consider $\bb_0$ a submanifold of the Euclidean space $\cc^n$ of dimension $m$ and codimension $k$. Its lift via $\pi$ is a submanifold $\mm_0$ invariant respect to vertical translations. 
Using notation \ref{notazione_esempi}, by \eqref{nbheis} and some trivial computation we have:
\begin{eqnarray*}
A'(X_j\hl,V)&=
   & \frac{1}{2}\sum_{\alpha=1}^k\bar g\left(X_j\hl,J\xi_{\alpha}\hl\right)\xi_{\alpha}\hl = \frac{1}{2}\sum_{\alpha=1}^k\bar g\left(X_j,J\xi_{\alpha}\right)\xi_{\alpha}\hl,
\end{eqnarray*}
where $J$ in the last term is the usual complex structure of $\cc^n$. This result is very similar to what we have for Hopf fibration.  It follows that
$$
\aaap=\aaa+\frac{1}{2}\sum_{\alpha=1}^{k}\left|J\xi_{\alpha}^{\top}\right|^2,
$$
then
\begin{equation}\label{aaa_h}
\aaa\leq\aaap\leq\aaa+\frac{k}{2}, 
\end{equation}
with $\aaap=\aaa$ if and only if for every $\alpha$ $J\xi_{\alpha}$ is normal to $\bb_0$, that is $\bb_0$ is a complex submanifold of $\cc^n$, while $\aaap=\aaa+\frac{k}{2}$ if and only if for every $\alpha$ $J\xi_{n+\alpha}$ is tangent to $\bb_0$, that is $\bb_0$ is CR-submanifold of $\cc^n$ of CR-dimension $m-k$. In the first case, in particular, $\bb_0$ is a minimal submanifold.  The classical Huisken's result \cite{H1} about evolution of convex hypersurfaces of the Euclidean space gives the following result for hypersurfaces of the Heisenberg group.

\begin{Proposition}\label{hei01}
Let $\mm_0$ an hypersurface of $\qq^n$. If $\mm_0$ is a cylinder with vertical axis, without boundary and its projection via $\pi$ is a convex hypersurface of $\mathbb R^{2n}$, then there is an unique solution of the mean curvature flow of $\mm_0$ invariant respect to vertical translations. Moreover this solution develops a singularity in finite time and converges to a vertical line. Then such an $\mm_0$ is diffeomorphic to a cylinder $\sss^{2n-1}\times\mathbb R$.
\end{Proposition}
\proof Such an $\mm_0$ is invariant respect to vertical translation, the fiber of $\pi$ are not closed so we can apply Theorem \ref{mainsub} in the sense of remark \ref{non_unico}. Let $\bb_0=\pi(\mm_0)$. By the main result of \cite{H1}, $\bb_0$ shrinks to a round point in finite time. The thesis follows lifting this result to $\mm_0$. \cvd

Using the main theorem of \cite{AB}, we have the following result for submanifolds of arbitrary codimension in the Heisenberg group.


\begin{Proposition}\label{hei02}
Let $\mm_0$ a cylinder with vertical axis of $\qq^n$ of dimension $m\geq 3$, without boundary and whose horizontal section is a closed submanifold. If $\mm_0$ has $H\neq 0$ everywhere and satisfies $\aaa\leq c \hhh$ with
$$
c\leq\left\{\begin{array}{ll}
\frac{4}{3(m-1)} & \quad\text{if }\quad 3\leq m \leq 5,\\
\frac{1}{m-2} & \quad\text{if }\quad m>5,
\end{array}\right.
$$
then the mean curvature flow of initial data $\mm_0$ has an unique $\mathbb R$-invariant solution and this solution converges in finite time to a vertical line. Hence such an $\mm_0$ is diffeomorphic to a cylinder $\sss^{m-1}\times \mathbb R$.
\end{Proposition}
\proof We have that $\bb_0=\pi(\mm_0)$ is a closed submanifold of $\mathbb R^{2n}$ of dimension $m-1$. By \eqref{aaa_h}, $\bb_0$ satisfies
$$
\aaa\leq\aaap\leq c\hhhp=c\hhh.
$$
The main result of \cite{AB} says that the evolution by mean curvature of $\bb_0$ shrinks to a point in finite time. We can apply Theorem \ref{mainsub} to the unique $\mathbb R$-invariant solution obtaining the convergence of $\mm_0$ to a fiber of $\pi$, that is to a vertical line of $\qq^n$.\cvd

Another interesting submersion is the one that arise with the tangent sphere bundle of a Riemannian manifold equipped with the Sasaki metric. For any Riemannian manifold $(\bbb,\bar g)$ let $T\bbb$ its tangent bundle and for any $r>0$ let $T^r\bbb=\left\{(p,u)\in T\bbb\left|\left|u\right|_{\bar g}=r\right.\right\}$ be the tangent sphere bundle of radius $r$. The natural projection
$$
\pi:(p,u)\in T^r\bbb\mapsto p\in\bbb
$$
is a submersion. In this special case, for any $X$ vector field on $\bbb$ we can define also a lift $X\tl\in\vc$ called \emph{tangent lift}: see \cite{KS} for an exhaustive description. The Sasaki metric is a natural metric $\bar g$ on $T\bbb$, restricted to $T^r\bbb$ has the following form:
\begin{equation}\label{metrica_sasaki}
\begin{array}{rcl}
\bar g_{(p,u)}(X\hl,Y\hl) & = & \displaystyle{\bar g_p(X,Y),}\\
\bar g_{(p,u)}(X\tl,Y\tl) & = & \displaystyle{\bar g_p(X,Y)-\frac{1}{r^2}\bar g_p(X,u)\bar g_p(Y,u),}\\
\bar g_{(p,u)}(X\hl,Y\tl) & = & 0,
\end{array}
\end{equation}
for any $X$ and $Y$ tangent to $\bbb$. With this metric the projection $\pi:T^r\bbb\to\bbb$ is a Riemannian submersion with fibers $\pi^{-1}(p)=T^r_p\bbb$, the sphere of radius $r$ tangent to $\bbb$ in $p$. The horizontal distribution of $\pi$ is generated by the horizontal lifts and the vertical distribution is generated by the tangential lifts. The group of isometries that we are considering acts only on the vectorial part as an isometry of  $T^r_p\bbb$ and is isomorphic to $O(n)$, where $n$ is the dimension of $\bbb$. Note that in this case the action of the group is not free in fact the orbits are not isometric to the group, but the quotient manifold $T^r\bbb /O(n)\equiv\bbb$ is a well defined manifold. The Levi-Civita connection of the Sasaki on metric on $T^r\bbb$ is 

\begin{Lemma}\label{nbsasaki}\cite{KS}
For any $X$ and $Y$ vector fields tangent to $\bbb$ we have:
\begin{itemize}
\item[1)] $\left(\bnb_{X\hl}Y\hl\right)_{(p,u)}=\left(\bnb_XY\right)\hl_{(p,u)}-\frac{1}{2}\left(\bar R_p(X,Y)u\right)\tl$,
\item[2)] $\left(\bnb_{X\hl}Y\tl\right)_{(p,u)}=\left(\bnb_XY\right)\tl_{(p,u)}+\frac{1}{2}\left(\bar R_p(u,Y)X\right)\hl$,
\item[3)] $\left(\bnb_{X\tl}Y\hl\right)_{(p,u)}=\frac{1}{2}\left(\bar R_p(u,X)Y\right)\hl$,
\item[4)] $\displaystyle{\left(\bnb_{X\tl}Y\tl\right)_{(p,u)}=-\frac{1}{r^2}\bar g_p(u,Y)X\tl}$,
\end{itemize}
where $\bar R$ is the Riemann curvature tensor of $\bbb$.
\end{Lemma}

The fibers are closed and last equation shows that they are also totally geodesic: $\hat A(X\tl,Y\tl)$ is the horizontal part of $\bnb_{X\tl}Y\tl$. From now on consider a submanifold $\bb_0$ of dimension $n$ and codimension $k$ and $\mm_0$ its $O(n+k)$-invariant lift to $T^r\bbb$. Since in this case we have a way to lift vector fields on $\bbb$ to vector fields tangent to the fibers, we modify notation \ref{notazione_esempi}. For any $p\in\bb_0$ and any $(p,u)\in\pi^{-1}\left\{p\right\}$, let $(X_1,\dots,X_n)$ an orthonormal basis tangent to $\bb_0$ in $p$, $(\xi_1,\dots,\xi_k)$ an orthonormal basis normal to $\bb_0$ in $p$ such that
$$
u=r\cos(\vartheta)X_1+r\sin(\vartheta)\xi_1,
$$
for some $\vartheta$. 
Let $Z=\sin(\vartheta)X_1-\cos(\vartheta)\xi_1$, then $(u,Z,X_2,\dots,X_n,\xi_2,\dots,\xi_k)$ is an orthogonal basis of $T_p\bb$. By \eqref{metrica_sasaki} we have that
$(X_1\hl,\dots,X_n\hl,Z\tl,X_2\tl,\dots,X_n\tl,\xi_2\tl,\dots,\xi_k\tl)$ is an orthonormal basis tangent to $\mm_0$ in $(p,u)$, while $(\xi_1\hl,\dots,\xi_k\hl)$ is an orthonormal basis normal to $\mm_0$ in $(p,u)$. As concrete example, consider $\bbb=\sss^{n+k}(c)$ the sphere of constant curvature $c>0$. By Lemma \ref{nbsasaki} we have
\begin{eqnarray*}
A'(X_i\hl,Z\tl)(p,u) & = & \frac{1}{2}\sum_{\alpha=1}^k\bar R_p(u,Z,X_i,\xi_{\alpha})\xi_{\alpha}\hl= -\frac{cr}{2}\delta_{i1}\xi_{1}\hl.
\end{eqnarray*}
Similarly
\begin{eqnarray*}
A'(X_i\hl,X_j\tl)(p,u) & = & -\frac{cr}{2}\sin(\vartheta)\delta_{ij}\xi_{1}\hl,\\
A'(X_i\hl,\xi_j\tl)(p,u) & = & \frac{cr}{2}\cos(\vartheta)\delta_{i1}\xi_{j}\hl.
\end{eqnarray*}

\noindent Then
\begin{eqnarray*}
\aaap (p,u) & = &\aaa(p)+\frac{c^2r^2}{2}\left(1+(n-1)\sin^2(\vartheta)+(k-1)\cos^2(\vartheta)\right)\\
 & = & \aaa(p)+\frac{c^2}{2}\left(r^2+(n-1)\left|u^{\perp}\right|^2+(k-1)\left|u^{\top}\right|^2\right),
\end{eqnarray*}
where $\perp$ (respectively $\top$) indicates the normal (respectively the tangent) component respect to $\bb_0$.
In particular we have
$$
\aaa(p)+\frac{c^2r^2}{2}min\left\{k,n\right\}\leq\aaap (p,u)\leq \aaa(p)+\frac{c^2r^2}{2}max\left\{k,n\right\}.
$$

Lifting the submanifolds of the sphere considered by Huisken \cite{H3} and Baker \cite{Ba} we have the following result as consequence of Theorem \ref{mainsub}.

\begin{Proposition}\label{sasaki01}
For any $r>0$, $n\geq 3$ and $k\geq 1$, let $\mm_0$ be a $2n+k-1$-dimensional $O(n+k)$-invariant submanifold of $T^r\sss^{n+k}(c)$. Suppose that $\mm_0$ satisfies the pinching condition
$$
\aaa<\frac{1}{n-1}\hhh +2c + \frac{c^2}{2}min\left\{k,n\right\},
$$
then the mean curvature flow with initial data $\mm_0$ converges in finite time to a fiber $\pi^{-1}(p)=T^r_p\sss^{n+k}(c)$ or the flow is defined for any time and converges to $\pi^{-1}(\sss^n(c))$ that is a minimal, but not totally geodesic, submanifold of $T^r\sss^{n+k}(c)$.
\end{Proposition}





\noindent {\bf Acknowledgments} 
The results of this paper are part of Giuseppe Pipoli's PhD thesis, written at the Department of Mathematics, University ``Sapienza'' of Rome. 
Giuseppe Pipoli was partially supported by PRIN07 ``Geometria Riemanniana e strutture differenziabili'' of MIUR (Italy) and Progetto universitario Univ. La Sapienza ``Geometria differenziale -- Applicazioni''.

\bigskip

\noindent Giuseppe Pipoli, Institut Fourier, Universit\'e Joseph Fourier (Grenoble I), UMR 5582, CNRS-UJF, 38402, Saint-Martin-d'H\`eres, France. E-mail: giuseppe.pipoli@ujf-grenoble.fr

\end{document}